\documentclass{amsart}

\usepackage[T1]{fontenc}
\newtheorem{defn}{Definition}[section]

\newtheorem{theorem}[defn]{Theorem}

\theoremstyle{remark}
\newtheorem{rmk}[defn]{Remark}

\usepackage{ulem}
\usepackage[utf8]{inputenc}
\usepackage[T1]{fontenc}
\usepackage[english]{babel}
\usepackage{amsfonts}
\usepackage{mathtools}
\usepackage{amsthm}
\usepackage{amssymb}
\usepackage{caption}
\usepackage{graphicx}
\usepackage{enumitem}
\usepackage{verbatim}
\DeclarePairedDelimiter{\norm}{\lVert}{\rVert}

\newcommand{\Rn}{\mathbb{R}^n}
\DeclareMathOperator{\R3}{\mathbb{R}^3}
\DeclareMathOperator{\S3}{\mathbb{S}^3}
\DeclareMathOperator{\H3}{\mathbb{H}^3}
\DeclareMathOperator\supp{supp}
\renewcommand{\phi}{\varphi}
\renewcommand{\epsilon}{\varepsilon}

\title{Diameter estimates for surfaces in conformally flat spaces}
\author{Marco Flaim}
\author{Christian Scharrer}
\address{Institut für Angewandte Mathematik, Universität Bonn, Bonn, Germany}
\email[Marco Flaim]{flaim@iam.uni-bonn.de}
\email[Christian Scharrer]{scharrer@iam.uni-bonn.de}

\date{}

\begin{document}
	\maketitle
	
	\begin{abstract}
		The aim of this paper is to give an upper bound for the intrinsic diameter of a surface with boundary immersed in a conformally flat three dimensional Riemannian manifold in terms of the integral of the mean curvature and of the length of its boundary. Of particular interest is the application of the inequality to minimal surfaces in the three-sphere and in the hyperbolic space. Here the result implies an a priori estimate for connected solutions of Plateau's problem, as well as a necessary condition on the boundary data for the existence of such solutions. The proof follows  a construction of Miura and uses a diameter bound for closed surfaces obtained by Topping and Wu--Zheng.
	\end{abstract}
	
	\section{Introduction}
	In 1993, Simon \cite{MR1243525} obtained an upper bound for the extrinsic diameter involving the $L^2$-norm of the mean curvature and the $L^1$-norm of the second fundamental form for closed surfaces and surfaces with boundary, respectively.
	
	In a similar spirit, Topping \cite{MR2410779} proved in 2008 an upper bound for the intrinsic diameter of closed manifolds immersed in $\Rn$. More precisely, he found that for any closed, connected, $m$-dimensional manifold $M$ immersed in $\Rn$,
	\begin{equation}\label{topping}
		d(M) \leq C(m) \int_M |H|^{m-1} d\mu,
	\end{equation}
	where $d(M)$ is the intrinsic diameter, $H$ the mean curvature, $\mu$ the volume form associated to $M$, and $C$ a constant depending only on $m$. The result uses the Michael--Simon inequality (see the original \cite{MR344978} from 1973 or a more recent proof in \cite{brendle2022proof}): for $m\geq 2$, for any compact, $m$-dimensional manifold immersed in $\Rn$, and $f\in C^1(M)$ non-negative, 
	\begin{equation}\label{eq:MS}
		\left(\int_M f^{\frac{m}{m-1}} d\mu\right)^{\frac{m-1}{m}} \leq C(m) \left( \int_M \left(|\nabla f| + f |H|\right) d\mu + \int_{\partial M} f d\sigma \right)
	\end{equation}
	with again $C$ only depending on $m$.
	
	This latter inequality was extended one year later (1974) by Hoffman--Spruck \cite{MR365424,MR397625} to the case where $M$ is immersed in a Riemannian manifold. In this case some additional assumptions on $M$ depending on the sectional curvatures and on the injectivity radius of the ambient are needed, see Theorem \ref{hoffspruck} for details. In particular, \eqref{eq:MS} remains valid for $M$ immersed into a Riemannian manifold with non-positive sectional curvatures and infinite injectivity radius.
	
	In 2010, Wu--Zheng \cite{MR2823054} proved the analogous result of \eqref{topping} in the general Riemannian setting, adapting the proof of Topping to the more general framework and using the Hoffman--Spruck inequality in place of Michael--Simon's one. They showed that \eqref{topping} remains valid under the assumptions of Hoffman--Spruck.
	Note that Hoffman--Spruck's inequality (as well as Wu--Zheng's diameter bound) implies in particular the non-existence of closed minimal submanifolds in a simply connected complete  space with non-positive sectional curvatures (such ambients have infinite injectivity radius by Cartan--Hadamard's theorem).
	
	A clever construction was done by Miura in 2020 \cite{miura2020diameter}, in order to get a diameter bound for surfaces \textit{with boundary} immersed in $\Rn$. The idea is to double the surface (we sometimes informally refer to it with 2$\Sigma$), glue the two copies together and then apply Topping's diameter bound to the resulting closed surface. More precisely, Miura showed that for any connected compact surface $\Sigma$ immersed in $\Rn$ it holds
	\begin{equation}\label{eq:Miura}
		d(\Sigma) \leq C(2)\left(2\int_{\Sigma} |H| d\mu + \pi \ell(\partial M)\right)
	\end{equation}
	where $\ell(\partial M)$ is the length of the boundary (which may be disconnected) and $C(2)$ is the constant of \eqref{topping}. This type of estimate implies non-existence of compact connected solutions to the Plateau problem for certain boundary data (say when when the boundary is made of two short, distant curves).
	
	The idea of doubling and gluing the surface strongly relies on the Euclidean structure, but one can observe that the same can be done when the ambient is $(\Rn,\overline{g})$, with $\overline{g}$ a general metric. In particular, we deal with the case when $\overline{g}=e^{2\phi}\delta$ with $\delta$ the Euclidean metric and $\phi\in C^\infty(\Rn)$, which is to say when $\overline{g}$ is a conformal change of $\delta$. This choice is motivated by the application of the result to $\S3$ and $\H3$, see Section \ref{applications}.
	In Theorem \ref{mainthm} we prove that \eqref{eq:Miura} remains valid for $\Sigma$ immersed in a conformally flat region of a complete Riemannian manifold provided the surface satisfies additional smallness assumptions (see \eqref{cond1.3}, \eqref{cond2.3}).

	\section{Proof of the main result}\label{mainsection}
	Let us first fix some notations and conventions used throughout the paper. When the ambient manifold is not $(\mathbb R^n,\delta)$, we denote it with $(N,\overline{g})$.
	We will mainly use $M\subset N$ for closed immersed manifolds and $\Sigma\subset N$ for immersed manifolds with boundary, and with $g$ we denote the metric induced by restriction. When we need to see how some quantities transform under a different metric in $\mathbb R^n$, we write subscript $\delta$ to say that it is taken with respect to the Euclidean metric $\delta$ (for example, $g_\delta$ is the metric induced on $M$ by $\delta$, while $g$ is the one induced by $\overline{g}$). Similarly, we will need to distinguish between $\mathrm D$, the Levi-Civita connection of $\delta$ and $\overline{\nabla}$, the one of $\overline{g}$. The volume form $d\mu$ is the one induced by $g$ on $M$ or $\Sigma$ and $|U|$ is the volume of some $U\subset M$ with respect to $d\mu$. The (vector-valued) second fundamental form is denoted with $A$ and the mean curvature with $H$, which is for us the trace of $A$. With $\eta$ we denote the unit normal field to the boundary $\partial\Sigma$, which is the vector field defined on the boundary tangent to the manifold, orthogonal to the boundary, pointing outward. When the immersed submanifold is a curve, we denote its extrinsic curvature with $\kappa$.
	For what concerns the intrinsic geometry, $K_N$ represents the sectional curvatures of $N$, and we write $K_N\leq K$ in the sense that at any point any sectional curvature is bounded from above by $K$. The injectivity radius of $N$ is denoted by $\overline{R}$, and $\overline{R}(A)=\inf_A \overline{R}$ for $A\subset N$.
	
	In this section we are going to state more precisely the results of the introduction and give the proof of the main theorem. We start giving the precise statement of Hoffman--Spruck's inequality.
	
	\begin{theorem}[Hoffman--Spruck]\label{hoffspruck}
		Let $M$ be a compact connected manifold  immersed in a Riemannian manifold $(N,\overline{g})$, assume $K_N \leq K\in\mathbb{R}$ and consider $f\in C^1(M)$ non-negative which is zero on the boundary. Define $\rho_0=\rho_0(\alpha,\supp(f))$ for $\alpha\in(0,1)$ by
		\[
		\rho_0=
		\begin{cases} K^{-1/2}\sin^{-1}\left[K^{1/2}(1-\alpha)^{-1/m}(\omega_m^{-1}|\supp f|)^{1/m}\right] & \text{for $K>0$} \\
			(1-\alpha)^{-1/m}(\omega_m^{-1}|\supp f|)^{1/m} & \text{for $K\leq 0$.}
		\end{cases}
		\]
		and assume the two following conditions:
		\begin{equation}\label{cond1}
			K(1-\alpha)^{-2/m}(\omega_m^{-1}|\supp f|)^{2/m} \leq 1 \tag{$\star$}
		\end{equation}
		\begin{equation}\label{cond2}
			2\rho_0\leq \overline{R}(M). \tag{$\star\star$}
		\end{equation}
		
		Then it holds
		\begin{equation*}
			\left(\int_M f^{\frac{m}{m-1}} d\mu\right)^{\frac{m-1}{m}} \leq c(m,\alpha) \int_M \left(|\nabla f| + f |H|\right) d\mu
		\end{equation*}
		for \[c(m,\alpha)=
		\begin{cases} \frac{\pi}{2} 2^{m-2}\alpha^{-1}(1-\alpha)^{-1/m}\frac{m}{m-1}\omega_m^{-1/m} & \text{for $K\geq0$} \\
			2^{m-2}\alpha^{-1}(1-\alpha)^{-1/m}\frac{m}{m-1}\omega_m^{-1/m} & \text{for $K<0$.}
		\end{cases}\]
	\end{theorem}
	
	\begin{rmk}\leavevmode
		\begin{itemize}
			\item In the case where $N$ is simply connected, complete, and has non-positive sectional curvature, by Cartan--Hadamard's theorem the injectivity radius is infinite, hence both conditions are trivially satisfied for any $\alpha\in(0,1)$.
			
			\item If there exist some $\alpha$ such that the condition is satisfied, we obviously aim to choose the one which minimises $c$. If there are no restrictions, the optimal constant is $c(m,\alpha)=c(m,\frac{m}{m+1})$.
		\end{itemize}
	\end{rmk}

	Wu--Zheng \cite{MR2823054}, after the work of Topping \cite{MR2410779}, proved that the latter implies the following diameter estimate.
	
	\begin{theorem}[Wu--Zheng]\label{wuzheng}
		Let $M$ be a closed connected manifold  immersed in a complete Riemannian manifold $(N,\overline{g})$ and assume $K_N \leq K\in\mathbb{R}$. For $\alpha\in(0,1)$, let 
		\[   \rho_0=\rho_0(\alpha,M)=
		\begin{cases} K^{-1/2}\sin^{-1}\left[K^{1/2}(1-\alpha)^{-1/m}(\omega_m^{-1}|M|)^{1/m}\right] & \text{for $K>0$} \\
			(1-\alpha)^{-1/m}(\omega_m^{-1}|M|)^{1/m} & \text{for $K\leq 0$}
		\end{cases}
		\]
		and assume that the following conditions hold
		\begin{equation}\label{cond1.2}
			K(1-\alpha)^{-2/m}(\omega_m^{-1}|M|)^{2/m} \leq 1 \tag{$\star'$}
		\end{equation}
		\begin{equation}\label{cond2.2}
			2\rho_0\leq \overline{R}(M).\tag{$\star\star'$}
		\end{equation}
		Then
		\begin{equation*}
			d(M)\leq C(m,\alpha)\int_M |H|^{m-1} d\mu
		\end{equation*}
		where $C(m,\alpha)$ are different constants from the previous theorem and for example we can take $C(2,\alpha)=\frac{576\pi}{\alpha^2(1-\alpha)}$.
	\end{theorem}
 
	We are now in position to state our result.
	\begin{theorem}\label{mainthm}
		Let $\Sigma$ be a compact connected surface immersed in an open subset $U$ of a complete Riemannian manifold $(N,\overline{g})$ such that $(U,\overline g)$ is isometric to $(V,e^{2\varphi}\delta)$ with $V\subset \mathbb R^n$ and $\phi\in C^\infty(V)$. Consider $K_N\leq K\in\mathbb R$, fix $\alpha\in(0,1)$, let
		\begin{equation*}
			\rho_0 =\rho_0( 2 \Sigma ) =
			\begin{cases} K^{-1/2}\sin^{-1}\left[K^{1/2}(1-\alpha)^{-1/2}(\omega_2^{-1}2|\Sigma|)^{1/2}\right] & \text{for $K>0$} \\
				(1-\alpha)^{-1/2}(\omega_2^{-1}2|\Sigma|)^{1/2} & \text{for $K\leq 0$}
			\end{cases}
		\end{equation*}
		and assume
		\begin{equation}\label{cond1.3}
			K ( 1 - \alpha)^{-1} \omega_2^{-1} 2 | \Sigma | < 1 \tag{$\star''$}
		\end{equation} 
		\begin{equation}\label{cond2.3}
			2\rho_0 < \overline{R}(\Sigma). \tag{$\star\star''$}
		\end{equation}
		Then
		\begin{equation*}
                d(\Sigma)\leq C(2,\alpha)\left[2\int_\Sigma |H|d\mu + \pi\ell(\partial\Sigma)\right]
            \end{equation*}%
        where $C(2,\alpha)$ is the best constant of Theorem \ref{wuzheng}.
	\end{theorem}
 
	\begin{rmk}\label{rmk:constant}
    The inequalities in conditions \eqref{cond1.3} and \eqref{cond2.3} need to be strict (in contrast with Theorem \ref{wuzheng}), if we want to take the best constant $C(2,\alpha)$ of Theorem \ref{wuzheng}. However, if we choose $C(2,\alpha)=\frac{576\pi}{\alpha^2(1-\alpha)}$, by continuity of $\alpha\mapsto C(2,\alpha)$, we can include the equality case in the conditions.
	\end{rmk}

    \begin{rmk}\label{rmk:hypothesis}
    Another possibility is to follow the approach of \cite{MR3293738}, where similar results to Theorem \ref{hoffspruck} and Theorem \ref{wuzheng} are proved for weighted Riemannian manifolds of the form $(\Rn,\delta,e^\psi d\mu_\delta)$. In this case, \eqref{cond1.2} and \eqref{cond2.2} are replaced by the condition 
        \begin{equation*}
            \int_\Sigma \langle x, \mathrm D\psi(x)\rangle_\delta e^{\psi(x)}d\mu_\delta(x) \geq 0.
        \end{equation*}
        However, we cannot adapt this result to the three-sphere because the condition is not satisfied by the weight relative to the stereographic projection. 
    \end{rmk}

	We now present the proof, starting by describing the construction of Miura \cite{miura2020diameter} for the Euclidean setting.
	The new part is to study how the quantities transform under the conformal change of the ambient metric and adapt the calculations.
	
	\begin{proof}
		According to our isometry assumption, we may treat $\Sigma$ as being immersed in $V\subset \mathbb R^n$. We would like to construct a closed surface from it, in order to apply Theorem \ref{wuzheng} and get some information on $\Sigma$. To do this, Miura doubled $\Sigma$ and glued the two copies in a suitable way. The idea is to pay as little as possible in terms of total mean curvature. 
		
		He first considered a "teardrop" curve \[z=(x,y):[0,b]\to\mathbb{R}^2\] immersed, closed ($z(0)=z(b)$), with $|z'|=1$ (here the norm is with respect to $\delta$ of $\mathbb{R}^2$) and such that $x'(0)=-x'(b)=1$, $y'(0)=y'(b)=0$. 
  
		Let also \[\gamma:[0,l]\to\Rn\] be a parametrization of $\partial\Sigma$ with $\gamma$ immersed, closed and $|\gamma'|_{\delta}=1$. On this curve consider a smooth orthonormal (with respect to $\delta$) moving frame $(e_1(t),e_2(t),e_3(t))$ such that $e_1(t)=\gamma'(t)$ and $e_2(t)=\eta(t)$ an outer unit normal vector on $\partial\Sigma$ (tangential to $\Sigma$).
		
		Now fix $\epsilon>0$ and define the surface given by the immersion \[F_\epsilon:[0,b]\times[0,l]\to\Rn\] with $F_{\epsilon}(s,t)=\gamma(t)+\epsilon(x(s)e_2(t)+y(s)e_3(t))$. We will denote the image of $F_\epsilon$ with $T_\epsilon$.
		
        It is clear that gluing $\Sigma$ with a copy of itself using this immersion, we get a $C^{1,1}$-regular closed immersed surface which we will call $M_\epsilon$. 
		Since $|T_\varepsilon|\to 0$ and $\rho_0(M_\epsilon)\to\rho_0(\Sigma)$ as $\varepsilon\to0$, we infer from \eqref{cond1.3} and \eqref{cond2.3} that 
			\[
			K(1-\alpha)^{-1}\omega_m^{-1}(2|\Sigma|+|T_\epsilon|) < 1 
			\]
			and
			\[
			2\rho(M_\varepsilon) < \overline{R}(M_\epsilon)
			\]
			for $\varepsilon$ small enough. Since these two inequalities are again strict, we can approximate $M_\varepsilon$ by smooth surfaces satisfying these inequalities and apply Theorem \ref{wuzheng} to deduce
		\begin{equation}\label{eq:diam-M_eps}
			d(M_\epsilon)\leq C(2,\alpha)\int_{M_\epsilon}|H|d\mu = C(2,\alpha)\left(2\int_{\Sigma}|H|d\mu + \int_{T_{\epsilon}}|H|d\mu \right).
		\end{equation}
		The idea is to use the trivial inequality $d(\Sigma)\leq d(M_\varepsilon)$, so that we can estimate the diameter of $\Sigma$ with the right-hand side of the previous inequality and ultimately send $\epsilon$ to zero. In particular we want to study the behaviour of 
		\[\int_{T_{\epsilon}}|H|d\mu.\]
  
		Our first step is to write $H$ (mean curvature vector with respect to $e^{2\varphi}\delta$) in terms of $H_\delta$ (with respect to $\delta$).
		By a standard formula, the Levi-Civita connection (recall that $\overline \nabla$ is the one of $\overline g = e^{2\varphi}\delta$, $\mathrm D$ of $\delta$) changes as \[\overline \nabla_XY = \mathrm D_X Y + X(\phi)Y + Y(\phi)X - \langle X,Y\rangle_\delta \mathrm D \phi\] so that the second fundamental form becomes 
		\begin{equation*}
			\begin{split}
				A(X,Y) &= \left(\overline\nabla_X Y\right)^\perp \\
				&= \left(\mathrm D_X Y + X(\phi)Y + Y(\phi)X - \langle X,Y\rangle_\delta \mathrm D \phi \right)^\perp \\
				&= \left( \mathrm D_X Y - \langle X,Y\rangle_\delta \mathrm D \phi \right)^\perp \\
				&=A_\delta(X,Y)-\langle X,Y\rangle_\delta\left(\mathrm D \phi \right)^\perp 
			\end{split}
		\end{equation*} 
		for any $X,Y\in T\Sigma$.
        Tracing, we see how the mean curvature vector changes:
		\begin{equation*}
			\begin{split}H &= g^{ij}A_{ij} \\
				&= e^{-2\phi} \delta^{ij} \left((A_\delta)_{ij} -\delta_{ij}\left(\mathrm D \phi \right)^\perp\right) \\
				&= e^{-2\phi}\left(H_\delta-2\left(\mathrm D\phi\right)^\perp\right)
			\end{split}
		\end{equation*}
  
		Since the norm changes as $|\cdot|=e^\phi |\cdot|_\delta$ and the volume form as $d\mu=e^{2\phi}d\mu_\delta$,
		\begin{align} \notag 
            \int_{T_{\epsilon}}|H|d\mu &= \int_{T_{\epsilon}} \left|H_\delta - 2\left(\mathrm D\phi\right)^\perp \right|_\delta e^{\phi}d\mu_\delta \\ \label{eq:integral}
            &\leq \int_{T_{\epsilon}}|H_\delta|_\delta e^{\phi}d\mu_\delta + \int_{T_{\epsilon}}\left|2\left(\mathrm D\phi\right)^\perp\right|_\delta e^{\phi}d\mu_\delta.    
		\end{align}
		Now, $\Sigma$ is compact and immersed in the open set $V$. Therefore we can choose a compact set $C$ such that $M_\varepsilon\subset C\subset V$ for $\varepsilon$ small. It follows
		\[\left|2\left(\mathrm D\phi\right)^\perp\right|_\delta e^{\phi} \leq 
        2\left|\mathrm D\phi\right|_\delta e^{\phi} \leq 
        2\norm{\mathrm D\phi}_{L^\infty(C)}e^{\norm{\phi}_{L^\infty(C)}}\]
		and, for $\epsilon\to0$, $|T_\epsilon|\to 0$ implies
		\[\int_{T_{\epsilon}}\left|2\left(\mathrm D\phi\right)^\perp\right|_\delta e^{\phi}d\mu \longrightarrow 0.\]
  
		For the first term of \eqref{eq:integral}, Miura computed the integral of the mean curvature of $T_\epsilon$, showing that $|H_\delta|_\delta d\mu_\delta=|\kappa_z(s)| dsdt + O(\epsilon)$, where $\kappa_z$ is the curvature of the curve $z$ in $(\mathbb{R}^2,\delta)$. Therefore, taking care of the conformal factor, we get:
		\begin{equation*}
			\int_{T_\epsilon}|H_\delta|_\delta e^{\phi}d\mu = \int_0^l \int_0^b |\kappa_z(s)|e^{\phi(F(s,t))}ds dt + O(\epsilon).
		\end{equation*}
		Moreover,
		\[e^{\phi(F_\epsilon(s,t))}=e^{\phi(F_\epsilon(0,t))}+O(\epsilon)=e^{\phi(\gamma(t))}+O(\epsilon).\] 
		Thus, abbreviating $\mathcal K(z) = \int_0^b|\kappa_z(s)|ds$,
		\begin{equation*}
			\int_{T_\epsilon}|H_\delta|_\delta e^{\phi}d\mu_\delta = \mathcal K(z) \int_0^l e^{\phi(\gamma(t))} dt + O(\epsilon) = \mathcal K(z) {\ell}(\gamma) + O(\epsilon)
		\end{equation*}
		where ${\ell}(\gamma)$ is the length with respect to $g$ because the curve $\gamma$ was parametrized with respect to the arclength of $\delta$ and
		\[e^{\phi(\gamma)}=e^{\phi(\gamma)}|\gamma'|_\delta = |\gamma'|_{g}.\]
  
		Altogether, \eqref{eq:diam-M_eps} becomes
		\[d(\Sigma) \leq C(2,\alpha)\left( 2\int_\Sigma |H| d\mu + \mathcal K(z)\ell(\partial\Sigma) \right).\]
		The conclusion follows since by \cite[Lemma 2.2]{miura2020diameter}, the teardrop curve $z$ can be chosen in a way that $\mathcal K(z)$ approaches $\pi$.
	\end{proof}
	
	\begin{rmk}
		Some efforts to extend Theorem \ref{wuzheng} to surfaces with boundary was already done by Paeng \cite{MR3183369} and Wu \cite{wu2022diameter}. While their theorems hold in the more general setting of Riemannian manifolds, without the assumption of being conformally flat, they require the surface $\Sigma$ to be geodesically convex. In particular, as opposed to our theorem, their results cannot be used to obtain a priori statements about solutions of Plateau's problem. Diameter bounds for Euclidean submanifolds with boundary were obtained in \cite{menne2017novel}.
	\end{rmk}

    \begin{rmk}
        Theorem \ref{mainthm} can be applied to locally conformally flat (LCF) Riemannian manifolds, provided that $\Sigma$ is so small that it can be covered by the domain of a conformally flat chart. The Weyl--Schouten theorem states that an $n$-dimensional manifold is LCF if and only if:
        \begin{itemize}
            \item when $n=3$, the Schouten tensor $S(X,Y)=\text{Ric}(X,Y)-\frac{R}{4}g(X,Y)$ satisfies $\overline\nabla_X S(Y,Z)=\overline\nabla_Y S(X,Z)$;
            \item when $n\geq 4$, the Weyl tensor vanishes.
        \end{itemize}
    \end{rmk}

	\section{Applications}\label{applications}
	In this section we want to explore the consequences of the obtained inequality. Assuming $\Sigma$ to be a minimal surface satisfying the assumptions of Theorem \ref{mainthm}, one gets
	\begin{equation}\label{sphere case}
	    d(\Sigma) \leq C(2,\alpha)\pi\ell(\partial\Sigma).
	\end{equation} 
	In particular, as already noticed by Miura, this gives a necessary condition for a given boundary $\Gamma$ to admit the existence of a connected solution for the Plateau's problem. In the Euclidean case, this implies the well known fact that, given two parallel circles, if we move them far enough apart, there does not exist a connected compact minimal surface spanning the two circles (one can imagine that the catenoid exists only when the circles are close to each other). When the ambient is a Riemannian manifold things are a bit more involved because of the additional conditions and the constant that we obtain is very large, but still we can derive qualitatively interesting conditions. In our examples we will use the result in the version of Remark \ref{rmk:constant}.
	
	\subsection{The hyperbolic space}
	The hyperbolic space $\H3$ is conformally flat, as one can see from the Poincaré disk-model 
	\[N=\{r\leq 1\}\subset\R3, \ \ \ \ \ \overline{g}=\frac{4}{(1-r^2)^2}\delta\]
	where $r=r(x,y,z)=\sqrt{x^2+y^2+z^2}$, or from the Poincaré half-plane model
	\[N=\{z>0\}\subset\R3, \ \ \ \ \ \overline{g}=\frac{1}{(z^2)}\delta.\]

	Moreover, the space has constant sectional curvatures $K=-1$ and unbounded injectivity radius. Hence the usual conditions \eqref{cond1.3} and \eqref{cond2.3} are satisfied by any $\Sigma$ and any $\alpha\in(0,1)$. As observed by Wu--Zheng, choosing the optimal $\alpha=2/3$, one gets $C(2,2/3)=3888\pi$. The inequality implies then that for any compact, connected minimal surface immersed in $\H3$,
	\[d(\Sigma)\leq 3888\pi\ell(\partial\Sigma).\]

	\subsection{The 3-sphere}
	Via the stereographic projection, we can identify the unit 3-sphere, after removing one point, with
	\[ V=\R3\setminus\{0\}, \ \ \ \ \ \overline{g}=\frac{4}{(1+r^2)^2}\delta. \]
 
	We have that $K=1$ and $\overline{R}=\pi$ at any point, hence condition \eqref{cond1.3} becomes
	\begin{equation*}
	    \alpha\leq 1-\frac{2|\Sigma|}{\pi}
	\end{equation*}
    while condition \eqref{cond2.3} $\rho_0\leq\pi/2$ is always satisfied.
	Therefore, when $|\Sigma|\leq\frac{\pi}{6}$ we can take $C(2,2/3)=3888\pi$ as above.
    In the plane, if a minimal surface does not satisfy \eqref{eq:Miura}, then it is disconnected or non-compact (or both). Analogously, in $\S3$, a minimal surface contradicting \eqref{sphere case} is disconnected or has area bigger than $\frac{\pi}{6}$. For solutions of Plateau's problem, one can exclude the case that has $|\Sigma|$ area bigger than $\frac{\pi}{6}$, provided there exists a competitor of area smaller than $\frac{\pi}{6}$.
	
	\bibliographystyle{alpha}
    \bibliography{bibliography}
\end{document}